\documentclass{amsart}

\usepackage{amscd}

\numberwithin{equation}{section}

\newtheorem{thm}[equation]{Theorem}
\newtheorem{lem}[equation]{Lemma}

\newtheorem{cor}[equation]{Corollary}
\newtheorem{prop}[equation]{Proposition}

\theoremstyle{definition}

\newtheorem{defn}[equation]{Definition}
\newtheorem{eg}[equation]{Example}

\theoremstyle{remark}

\newtheorem{rmk}[equation]{Remark}

\DeclareMathOperator{\minp}{min.poly.}
\DeclareMathOperator{\gminp}{ch.poly.}
\DeclareMathOperator{\chpoly}{ch.poly.}

\newcommand{\Id}{\mathrm{Id}}

\newcommand{\tbtmat}[4]{\begin{pmatrix} #1 & #2 \\ #3 & #4 \end{pmatrix}}

\renewcommand{\H}{\mathbb{H}}
\newcommand{\R}{\mathbb{R}}
\newcommand{\C}{\mathbb{C}}
\newcommand{\Q}{\mathbb{Q}}
\newcommand{\Z}{\mathbb{Z}}
\newcommand{\F}{\mathbb{F}}

\newcommand{\ra}{\rightarrow}

\newcommand{\Fbar}{\overline{F}}

\newcommand{\la}{\lambda}

\DeclareMathOperator{\End}{End} \DeclareMathOperator{\adj}{adj}
\DeclareMathOperator{\sign}{sign}
\DeclareMathOperator{\cdet}{Cdet}
\DeclareMathOperator{\tr}{tr}
\DeclareMathOperator{\Pf}{Pf}
\DeclareMathOperator{\Skew}{Skew}
\DeclareMathOperator{\chr}{char}
\renewcommand{\det}{\mathrm{det}}

\newcommand{\EndF}{\End_F}
\newcommand{\EndK}{\End_K}

\newcommand{\e}{\varepsilon}
\newcommand{\s}{\sigma}

\newcommand{\ot}{\otimes}

\newcommand{\mysubsection}[1]{%
\medskip\noindent\refstepcounter{equation}\textbf{\theequation.\ #1.}}

\newcommand{\borelless}{%
\medskip\noindent\refstepcounter{equation}\textbf{\theequation.}~}

\begin{document}

\nocite{Axler} \nocite{Wstr} \nocite{Garlic}

\title[the characteristic polynomial]{The characteristic polynomial and determinant are not ad hoc constructions}

\author{R.~Skip Garibaldi}
\address{Dept.~of Mathematics and Computer Science, Emory University, Atlanta, GA 30322, USA}
\email{skip@member.ams.org}
\urladdr{http://www.mathcs.emory.edu/{\textasciitilde}skip}


\date{\today}

\maketitle

Most people are first introduced to the characteristic polynomial
and determinant of a matrix in a linear algebra course as
undergraduates.  The determinant is
usually defined as an alternating sum of products of entries of
the matrix (as in Jacobi \cite[\S4]{Jacobi}, \cite[\S23]{Weber:I},
and \cite[\S 6.1]{Bretscher}) or as the unique map $M_n(F) \ra F$
which is multilinear and alternating in the columns and which is 1
on the identity matrix (as in Weierstrass \cite{Wstr} and the
books by Hungerford \cite{Hung}; Lang \cite{Lang}; and Dummit and Foote).  
As a student, I thought these definitions were at best magical and at worst ad hoc.  
Where did the determinant come from?  
This paper gives definitions which I hope the reader will find more natural.

\smallskip

Now, the determinant of a linear transformation on $\R^n$ is a
natural enough object: its absolute value gives the factor by
which the transformation enlarges volumes, and its sign says
whether or not the map preserves orientation.  These properties
imply Weierstrass' axioms, see e.g.~\cite{Hannah} or \cite[\S5]{Lax}.

Another good definition of the determinant --- not so common at the undergraduate level --- is in terms of the $n$-th exterior power of $\wedge^n F^n$ as in \cite[Ch.~III, \S
6]{Bbk:alg}. This also leads to the Weierstrass axioms.

But even these two ``good'' definitions have a taint of being special to matrices.  (The first is even limited to matrices with real entries.)  After all, analogues of the determinant are known for the quaternions, the octonions, finite-degree field extensions...  It is not clear how to adapt the two good definitions to handle these algebras.  As mathematicians, we should demand a definition that works simultaneously in all cases.  Here we give a definition of the characteristic polynomial in \S\ref{chpoly.sec} which works for all of these cases; the constant term of this characteristic polynomial gives an analogue of the determinant.  For $n$-by-$n$ matrices, we derive Jacobi's alternating sum formula for the determinant, see \ref{alt.sum}.2.  We also recover the known ad hoc formulas for the determinant for quaternions (in \ref{quats.eg}) and finite-degree field extensions (in \ref{fields.eg}).  Moreover, the product formula $\det(aa') = \det(a) \det(a')$ always holds, see \S\ref{prod.sec}.

(One normally begins with a definition for the determinant and then defines the characteristic polynomial of a matrix $a$ as $\det(x - a)$.  We work in the opposite direction here.)

\smallskip

The
philosophy is the following.  Consider the lines in $\R^2$ given by the equations
  \[
   ax + by = c \quad \text{and} \quad a' x + b' y = c'.
  \]
If the coefficients $a$, $b$, $c$, $a'$, $b'$, $c'$ are specific real numbers, the lines might be parallel or the same (degenerate case), but ``typically'' they intersect at exactly one point.  If we treat the coefficients as independent indeterminates, we say that the lines are {\em generic}.  Such lines intersect at the point
  \[
  (x, y) = \left( \frac{b' c - c' b}{a b' - a' b}, \frac{a c' - a' c}{a b' - a' b}\right).
  \]

Similarly, a typical $n$-by-$n$ matrix has
$n$ distinct eigenvalues.\footnote{If you randomly choose an $n$-by-$n$ real matrix, the probability that you pick one with distinct eigenvalues is $100\%$.  This is because the other real matrices are a set with Lebesgue measure 0.  

Topologically, amongst the $n$-by-$n$ matrices, those with distinct
eigenvalues form a dense open subset.    Over $\R$ or $\mathbb{C}$
this is true in the usual topology; it is also true over any
infinite field in the Zariski topology \cite[\S{I.9}]{Jac:PI}.} For such a
matrix, the traditional characteristic polynomial is just the
minimal polynomial.  To define the characteristic polynomial of a particular matrix $a$, we first find the minimal polynomial $\minp_\gamma(x)$ of a generic matrix $\gamma$.  Plugging in specific values for the indeterminates in $\gamma$, we get a polynomial whose only indeterminate is $x$, and this is the characteristic polynomial of $a$.
This method of defining the
characteristic polynomial works for all finite-dimensional
$F$-algebras, and the determinant is (up to a sign) the constant
term of the characteristic polynomial.

\medskip

The core of the idea --- looking at the minimal polynomial of a generic element --- goes back to the late 1800s, see \cite[p.~241]{Study} and \cite[p.~301]{Scheffers}.  
However, all treatments that I have
found do not develop the properties of the general characteristic
polynomial (as in \cite{Dickson} and \cite{Deuring}), or they make
use of known properties of the characteristic polynomial and
determinant for matrices in studying the general characteristic
polynomial (as in \cite[\S I.13]{Albert}, \cite{Jac:J},
\cite{Jac:gn}, and \cite[\S5.18]{Jac:tor}).  We use only
elementary properties of matrices from the very nice paper
\cite{Axler:dwd}.

\smallskip

Readers with an algebraic background may argue that one can obtain
the characteristic polynomial of an $n$-by-$n$ matrix over an
arbitrary field $F$ by applying the structure theory for finitely-generated torsion modules over a PID (as is done in
\cite[\S6.7]{Herstein}).  But once one is using that much algebra,
the contents of this paper are not so far away, and the results
here are much stronger.

\section{The characteristic polynomial} \label{chpoly.sec}

In this section, we define the characteristic polynomial of an element $a$ in a finite-dimensional $F$-algebra $A$ (\ref{poly.def}) and give some of its basic properties (\ref {Jac.prop}).

\begin{defn} An {\em $F$-algebra} is a ring $A$ with a multiplicative identity $1 \ne 0$ such that $A$ is
an $F$-vector space and $\alpha(ab) = (\alpha a)b = a
(\alpha b)$ for all $\alpha \in F$ and $a, b \in A$. (Alternately,
$A$ is a ring with identity such that there is a monomorphism $F
\ra Z(A)$ which maps the identity in $F$ to the identity in $A$.)
All algebras that we consider will be {\em finite-dimensional} as
vector spaces over $F$.
\end{defn}

Principal examples of $F$-algebras are the $n$-by-$n$ matrices
$M_n(F)$ and Hamilton's quaternions, which form an $\R$-algebra (see \ref{quats.eg}).

Note that the definition gives a copy of $F$ inside the center of $A$, but there is no requirement that $F$ is the entire center of $A$; for
example, $A$ may be taken to be a finite-degree field extension of $F$.

We do require that $A$ is associative, e.g., $A$ cannot be the octonions.  However, the definition of the characteristic polynomial given below also works for the much broader
class of strictly power-associative $F$-algebras, which includes the octonions and
Jordan algebras.  Many of the properties of the characteristic polynomial proved here can also be proved in that more general setting, see \cite[Ch.~VI]{Jac:J}.

The only tool we need that is not always discussed in a first-year
graduate algebra course is the tensor product $\otimes$, for which
we refer the reader to any good graduate algebra text.
Heuristically, it allows one to make precise the notion of
enlarging our base field $F$: if $K$ is any extension field of
$F$, then $M_n(F) \otimes_F K$ is isomorphic to $M_n(K)$ as
$K$-algebras.  For $A$ an $F$-vector space, $\dim_F A = \dim_K (A \ot_F K)$.

Let $a_1, \ldots, a_m$ be an $F$-basis for $A$.  Let $R = F[t_1,
\ldots,  t_m]$ for $t_1, \ldots, t_m$ (commuting) indeterminates, and let $K$
be the quotient field of $R$.  We call $\gamma = \sum_i t_i a_i\in A \otimes_F K$ a {\em generic element}.  The $K$-span of $1,
\gamma, \gamma^2, \ldots$ is a subspace of $A \otimes K$, so it
must be finite-dimensional over $K$.  Hence there is a nonzero
monic polynomial $\minp_{\gamma/K}$ in $K[x]$ of smallest degree
such that $\minp_{\gamma/K}(\gamma) = 0$, called the {\em minimal
polynomial} for $\gamma$ over $K$.

Note that this polynomial is unique: If $f(x)$ and $g(x)$ are monic polynomials of minimal degree such that $f(\gamma) = g(\gamma) = 0$, then $h(x) = f(x) - g(x)$ is a polynomial of smaller degree such that $h(\gamma) = 0$.  This contradicts the minimality of $f$ and $g$ unless $h(x) = 0$.

\begin{lem} \label{poly.lem}
The minimal polynomial $\minp_{\gamma/K}$ is in $R[x]$, not just
$K[x]$.
\end{lem}

\begin{proof} 
Consider the $R$-submodules $A_j$ of $A \otimes R$ generated by $\{ 1,
\gamma, \gamma^2, \ldots, \gamma^j \}$.  They form an ascending
chain $A_1 \subseteq A_2 \subseteq \cdots$.  Since $R$ is
noetherian (Hilbert's Basis Theorem) and $A \otimes R$ is a
finitely-generated $R$-module, this chain must stabilize.  That
is, $\gamma_{j+1}$ is in $A_j$ for some $j$, so $\gamma$ satisfies
a monic polynomial $f$ in $R[x]$.\footnote{A more direct argument
would be: $R[\gamma]$ is an $R$-submodule of $A \otimes R$ and $A \otimes R$ is a finitely-generated $R$-module.  Hence $\gamma$ is
integral over $R$ \cite[Thm.~VIII.5.3]{Hung}.  Unfortunately, the
typical proof of this implication invokes determinants, so we use
instead that $R$ is noetherian.} Since $\minp_{\gamma/K}$ divides
$f$ in $K[x]$ and both are monic, $\minp_{\gamma/K}$ lies in $R[x]$ by
Gauss' Lemma.
\end{proof}

\begin{defn} \label{poly.def}
Write $a \in A$ with respect to the basis $a_1, \ldots, a_m$ above
as $a = \sum_i \alpha_i a_i$ for $\alpha_i \in F$.  The substitution $t_i \mapsto
\alpha_i$ defines a map $R[x] \ra F[x]$.  We call the image of
$\minp_{\gamma/K}$ in $F[x]$ the {\em characteristic polynomial of
$a$} and denote it by
   $\gminp_{a, A/F}$ or simply $\gminp_a$.
\end{defn}

\begin{rmk} \label{deg.rmk}
It is immediate from the definition that
   \[
   \deg(\gminp_{a,A/F}) \le \dim_F A
   \]
and that 
$ \deg(\gminp_{a,A/F})$
is the same for all $a  \in A$.
\end{rmk}

\begin{eg}[Upper-triangular matrices] \label{ut.eg}
Let $A$ be the algebra of $n$-by-$n$ upper triangular matrices
over $F$.  Write $E_{ij}$ for the matrix whose only nonzero entry
is a 1 in the $(i,j)$-position.  Fix a basis $a_1$, $\ldots$, $a_m$
for $A$ over $F$ consisting of $E_{ij}$'s with $a_i = E_{ii}$
 for $1 \le i \le n$.  Let $\gamma$ be the generic
element defined above.

Let $I_n$ denote the $n$-by-$n$ identity matrix.  For each $1 \le
i \le n$, the matrix $\gamma - t_i I_n$ has $n - 1$ pivot columns
--- equivalently, $n - 1$ leading 1's ---
in its row-reduced form, hence it has a nonzero kernel. That is,
$\gamma$ has an eigenvector in $K^n$ with eigenvalue $t_i$. Since
the $t_i$ are distinct elements of $K$, these eigenvectors form a
basis for $K^n$ and $\gamma$ is similar in $M_n(K)$ to the
diagonal matrix with diagonal entries $t_1$, $t_2$, $\ldots$,
$t_n$.  The minimal polynomial of $\gamma$ is $\prod_{i=1}^n (x -
t_i)$, since similar matrices have the same minimal
polynomials.\footnote{This argument may appear to be excessively
long. It is included here to illustrate that we are not making use
of determinants.}

By substitution, an upper triangular matrix $b \in A$ has
characteristic polynomial $\prod_{i = 1}^n (x - b_{ii})$.
\end{eg}

\begin{prop}
The characteristic polynomial $\gminp_{a,A/F}$ depends only on
$a$, $A$, and $F$ (and not on the choice of basis for $A$).
\end{prop}

\begin{proof}
Suppose that we have another $F$-basis $b_1, \ldots, b_m$ of $A$
with a corresponding generic element $\e = \sum_i t_i b_i$.  We
may write $b_i = \sum_j g_{ij} a_j$ for $g$ an invertible
$m$-by-$m$ matrix in $M_m(F)$.  Let $f \!: R \ra R$ be the $F$-algebra automorphism
defined by
  \[
  f(t_j) = \sum_i t_i g_{ij}.
  \]

Write $a \in A$ in terms of both bases as
  \begin{equation} \label{inv.1}
  a = \sum_j \alpha_j a_j = \sum_i \beta_i b_i.
  \end{equation}
We have a diagram
  \[
  \begin{CD}
  R[x] @>{t_j \mapsto \alpha_j}>> F[x] \\
  @V{f}VV @| \\
  R[x] @>{t_i \mapsto \beta_i}>> F[x]
  \end{CD}
  \]
with horizontal arrows the substitution maps.  Equation
\eqref{inv.1} gives that $\alpha_j = \sum_i \beta_i g_{ij}$ for
all $j$, hence the diagram commutes.

If we begin with $\minp_{\gamma/K}$ in the upper left,
substitution gives $\gminp_a$ computed with respect to the basis
$a_1, \ldots, a_m$ in $F[x]$.  On the other hand, $f$ extends in
an obvious way to an automorphism of $A \otimes R$ such that
  \[
  f(\gamma) = \sum_j f(t_j) a_j = \sum_j \Big(\sum_i t_i g_{ij}\Big) a_j =
  \sum_i t_i \Big(\sum_j g_{ij} a_j\Big) = \e.
  \]
Hence $f(\minp_{\gamma/K}) = \minp_{\e/K}$.  The image of this in
$F[x]$ is $\gminp_a$ computed with respect to the basis $b_1,
\ldots, b_m$.  The commutativity of the diagram gives the claim.
\end{proof}

\begin{lem} \label{scalar.ext}
Let $E$ be a field containing $F$ and fix $a \in A$.  The minimal polynomial and characteristic polynomial of $a$ is the same over $F$ and over $E$.
\end{lem}

The following precise statement of the lemma and its proof are technical and are best skipped by the casual reader.

\begin{proof}
More precisely we want to prove:
  \[
  \minp_{a/F} = \minp_{(a \otimes 1)/E}
  \]
and
  \[
  \gminp_{a,A/F} = \gminp_{(a \otimes 1), (A \otimes E)/E}.
  \]

Since $F$ is a field, $a$ generates a free submodule $B$ of $A$
with basis $1$, $a$, $a^2$, $\ldots$, $a^{d-1}$ for $d$ the degree
of $\minp_{a/F}$.  Then $B \otimes E$ is a free $E$-module with
the same basis, hence the degree of $\minp_{(a \otimes 1)/E}$ is
$\ge d$.  Since this polynomial divides $\minp_{a/F}$, they are
the same.\footnote{This paragraph may be replaced by the sentence:
Field extensions are ``faithfully flat''.  See \cite[pp.~45,
46]{AM} for a definition.}

The $F$-basis $a_1$, $\ldots$, $a_m$ of $A$ gives an $E$-basis
$a_1 \otimes 1$, $\ldots$, $a_m \otimes 1$ of $A \otimes E$, and
the generic element constructed from this $E$-basis is $\gamma
\otimes 1 \in E(t_1, \ldots, t_m)$.  Since the minimal polynomials of
$\gamma/K$ and $(\gamma \otimes 1)/E(t_1, \ldots, t_m)$ are the same by
the preceding paragraph, we get
   \[
   \gminp_{a, A/F} = \gminp_{(a \otimes 1), (A \otimes E)/E}
   \]
by substitution.
\end{proof}

In general, we write
  \begin{multline}  \label{char.poly}
  \gminp_a(x) = x^n - c_1(a) x^{n-1} + \cdots \\ \cdots + (-1)^{n-1} c_{n-1}(a) x
  + (-1)^n c_n(a).
  \end{multline}
The elements $c_1(a)$ and $c_n(a)$ play the roles of the trace and
determinant of $a$.

\begin{prop} \label{Jac.prop}
Let $A$ be a finite-dimensional $F$-algebra.  Then:
   \begin{enumerate}
   \item {\rm (Cayley-Hamilton)} $\gminp_a(a) = 0$ for all $a \in A$.
   \item If $f \!: A \ra A$ is a ring automorphism or anti-automorphism
   of $A$ which restricts to be an automorphism of $F$, then $f(c_i(a)) = c_i(f(a))$ for all $a \in A$.
   \item $c_i(\alpha a) = \alpha^i c_i(a)$ for $\alpha \in F$ and
   $a \in A$.
   \item $c_1 \!: A \ra F$ is $F$-linear.
   \item If $B$ is a subalgebra of $A$ and $b$ is in $B$, then
   $\gminp_{b,B/F}$ divides $\gminp_{b,A/F}$ in $F[x]$.
   \end{enumerate}
\end{prop}

We will observe in \ref{same.thm} below that our notion of
characteristic polynomial on $M_n(F)$ is the same as the usual
one. Then Prop.~\ref{Jac.prop} contains many results that one
typically proves in a linear algebra course. For example, (2)
gives that $\det(a^t) = \det(a)$ for $a \in M_n(F)$ and that
similar matrices have the same characteristic polynomial.

\begin{proof}
(1): Write $a = \sum_i \alpha_i a_i$.  Then $\gminp_a(a)$ is the
image of $\minp_{\gamma/K}(\gamma) \in A \otimes R$ under the
substitution map $A \otimes R \ra A \otimes F = A$ given by $t_i
\mapsto \alpha_i$.  Since $\minp_{\gamma/K}(\gamma) = 0$ in $A
\otimes R$, we have $\gminp_a(a) = 0$ in $A$.

\smallskip

(2): Suppose first that $f$ is an automorphism of $A$.  The
diagram
 \[
 \begin{CD}
 A \otimes R @>{t_i \mapsto \alpha_i}>> F[x] \\
 @V{f \otimes \Id_R}VV @VV{f}V \\
 A \otimes R @>{t_i \mapsto f(\alpha_i)}>> F[x]
 \end{CD}
 \]
commutes.  If we begin with $\minp_{\gamma/K}$ in the upper left,
we obtain $\gminp_a$ in the upper right $F[x]$, and then
$f(\gminp_a)$ in the lower right $F[x]$.  On the other hand, we
obtain $f(\minp_{\gamma/K}) = \minp_{f(\gamma)/K}$ in the lower left and then
$\gminp_{f(a)}$ in the lower right.  That is, we have the desired
equality $\gminp_{f(a)} = f(\gminp_a)$.

The same argument works in the case where $f$ is a ring
anti-automorphism, except that in the diagram we must replace $A
\otimes R$ in the lower left corner with $A^{\mathrm{op}} \otimes
R$ where $A^{\mathrm{op}}$ denotes the $F$-algebra $A$ with multiplication reversed.

\smallskip


(3) and (4): Suppose first that $\alpha$ is not 0.  If we write
the minimal polynomial of $\gamma/K$ as $\sum_{i=0}^n c_i x^i$ for
$c_i \in R$, then the minimal polynomial of $\alpha \gamma /K$ is
$\sum_{i=0}^n c_i \alpha^{n-i} x^i$.  Thus
  \[
  \alpha^n \minp_{\gamma/K}(x) = \minp_{\alpha \gamma/K}(\alpha x)
  \quad \text{in $R[x]$,}
  \]
and
  \[
  \alpha^n \gminp_a(x) = \gminp_{\alpha a}(\alpha x) \quad
  \text{in $F[x]$.}
  \]
Then we have $\alpha^n c_i(a) = \alpha^{n-i} c_i(\alpha a)$ for $0
\le i \le n$, which proves (3) for $\alpha \ne 0$. Since $c_i \!:
A \ra F$ is given by a polynomial in the coordinates of $a \in A$
with respect to some basis $a_1, \ldots, a_m$, this polynomial is
homogeneous of degree $i$.  This gives (4), as well as (3) for
$\alpha = 0$.

\smallskip

(5): Fix a basis $b_1$, $\ldots$, $b_r$ of $B$ and extend it to a
basis $a_1$, $\ldots$, $a_m$ of $A$ with $a_i = b_i$ for $1 \le i
\le r$.  By analogy, set $S = F[t_1, \ldots, t_r]$, let $L$ denote
the quotient field of $S$, and let $\e$ be the generic element
$\sum_{i=1}^r t_i b_i$ in $A \otimes S$.

Write $b = \sum_i \beta_i b_i$ for $\beta_i \in F$.  Let $\phi \!:
R \ra S$ be given by sending $t_j \mapsto 0$ for $r < j \le m$.
The image of $\minp_{\gamma/K}$ under the composition
   \[
   \begin{CD}
  R[x] @>\phi>> S[x] @>{t_i \mapsto \beta_i}>> F[x]
  \end{CD}
  \]
is $\gminp_{b,A/F}$.  Similarly, the image of $\minp_{\e/L} \in S[x]$ is $\chpoly_{b,B/F} \in F[x]$.

The homomorphism $\phi$ extends naturally to a map $A \ot R \ra A \ot S$ such that $\phi(\gamma) = \varepsilon$.  We have
  \[
  0 = \phi(0) = \phi(\minp_\gamma(\gamma)) = \phi(\minp_\gamma)\phi(\gamma),
  \]
which is the polynomial $\phi(\minp_\gamma) \in S[x]$ evaluated at $\phi(\gamma) = \varepsilon$.  Consequently, $\minp_{\e/L}$ divides $\phi(\minp_{\gamma/K})$ in $L[x]$.  Since $S$ is a UFD, $\minp_{\e/L}$ divides $\phi(\minp_{\gamma/K})$ in $S[x]$.  
Consequently, the image $\chpoly_{b,B/F}$ of $\minp_{\e,L}$ in $F[x]$ divides the image $\chpoly_{b,A/F}$ of $\phi(\minp_{\gamma/K})$.
\end{proof}


\section{Matrices}

In this section, we observe that the characteristic polynomial as
defined above agrees with the usual linear algebra notion of
characteristic polynomial in the case where $A = M_n(F)$.

Everyone knows the next lemma, but maybe not the clean proof:

\begin{lem} \label{minp.lem}
Let $T$ be a linear transformation of an $F$-vector space of dimension $n \ge 1$.  Then $T$ satisfies a nonzero polynomial of degree $\le n$.
\end{lem}

\begin{proof}  We sketch the nice proof from \cite{Burrow}.  Let $v$ be a nonzero vector in the vector space $V$.  The $n + 1$ vectors
   \[
   v, T(v), T^2(v), \ldots, T^n (v)
   \]
must be linearly dependent, so there is a polynomial $g(x) \in F[x]$ of degree $\le n$ such that $g(T)v = 0$.

Set $U := \ker g(T)$.  The linear transformations $T$ and $g(T)$ commute, so $T(U) \subseteq U$ and $T$ induces a linear transformation $T_{V/U}$ on $V/U$.  By induction, $T|_U$ satisfies a polynomial $m_U(x)$ of degree $\le \dim U$ and $T_{V/U}$ satisfies a polynomial $m_{V/U}(x)$ of degree $\le \dim V/U$.  Then $m_{V/U}(T) V \subseteq U$ and $T$ satisfies the polynomial $m_U(x)\, m_{V/U}(x)$.  Moreover,
   \[
   \deg (m_U \cdot  m_{V/U}) = \deg(m_U) + \deg(m_{V/U}) \le \dim U + \dim (V/U) = \dim V.\qedhere
  \]
\end{proof}

\begin{cor} \label{deg.lem}
The characteristic polynomial (in our sense) of a matrix in
$M_n(F)$ has degree $\le n$.
\end{cor}

\begin{proof}
Lemma \ref{minp.lem}  applies in particular to
$\gamma \in M_n(K)$ as a linear transformation of $K^n$.  Substitution gives the corollary.
\end{proof}

Fix an algebraic closure $\Fbar$ of $F$.  For $a \in M_n(F)$, we
call $\la \in \Fbar$ an {\em eigenvalue} of $a$ if the kernel of
$(\la I_n - a)$ is nonzero. Let $U_\la$ denote the corresponding
generalized eigenspace, i.e., the set of vectors $v \in \Fbar^n$
lying in the kernel of $(\la I_n - a)^r$ for some natural number
$r$.  The {\em multiplicity} $m(\la)$ of an eigenvalue $\la$ is
$\dim_F U_\la$.

\begin{thm} \label{same.thm}
For $a \in M_n(F)$, the characteristic polynomial (as defined in
\ref{poly.def}) of $a$ factors in $\Fbar[x]$ as
   \begin{equation} \label{same.eq}
   (x - \la_1)^{m(\la_1)} (x - \la_2)^{m(\la_2)} \cdots (x -
   \la_k)^{m(\la_k)},
   \end{equation}
where $\la_1, \ldots, \la_k$ are the distinct eigenvalues of $a$.
\end{thm}

\begin{proof}
Let $B$ be the subalgebra of $M_n(F)$ consisting of upper
triangular matrices and let $a'$ be in $B$.  The characteristic polynomial $\chpoly_{a',B/F}$ of $a'$ as an element of $B$ was computed in Example \ref{ut.eg}; it is of the form \eqref{same.eq}.  By \ref{Jac.prop}.5, $\chpoly_{a',B/F}$ divides the characteristic polynomial $\chpoly_{a',M_n(F)/F}$ of $a'$ as an element of $M_n(F)$.  Since both polynomials are monic and have degree $n$ (by \ref{ut.eg} and \ref{deg.lem}), the theorem holds for upper triangular matrices.

Since the characteristic polynomial of the given matrix $a$ is
unchanged under scalar extension, we may assume that $F$ is algebraically closed, i.e., where
$F = \Fbar$.  Here we need one somewhat sophisticated result from
linear algebra: Since $F$ is algebraically closed, $a$ is similar
to an upper triangular matrix $a'$ \cite[Thm.~6.2]{Axler:dwd}. But
the theorem holds for $a'$ by the preceding paragraph. Since similarity changes neither the characteristic
polynomial (\ref{Jac.prop}.2) nor the eigenvalues, the theorem
holds for $a$.
\end{proof}

In \cite{Axler:dwd}, Axler develops many of the typical properties
of matrices (e.g., the existence of eigenvalues and the
decomposition with respect to generalized eigenspaces) over an
algebraically closed field without use of the determinant.  For
example, in \S5 of that paper he defines the characteristic
polynomial to be exactly the product displayed in the theorem.
Logically, one could insert the contents of this paper at that
point in his. 

\begin{cor} \label{divides.1} For $a \in M_n(F)$, 
the minimal polynomial $\minp_a$ and the characteristic polynomial $\gminp_a$ have the same irreducible factors in $F[x]$.
\end{cor}

\begin{proof}
Irreducible polynomials in $F[x]$ are determined (up to a scalar factor) by their roots in an algebraic closure $\Fbar$.  Thus we may assume that $F$ is algebraically closed.  

As in the proof of Th.~\ref{same.thm}, $a$ is similar to an upper triangular matrix since $F$ is algebraically closed.  Since conjugation changes neither the characteristic nor the minimal polynomial, we may also assume that $a$ is upper triangular.

By Th.~\ref{same.thm}, every irreducible factor of the characteristic polynomial is of the form $(x - \la_i)$ where $\la_i$ is a diagonal entry in $a$, say $\la_i = a_{ii}$.  The $(i, i)$-entry of $\minp_a(a) = 0$ is 0, but it is also $\minp_a(\la_i)$.  Therefore, $(x - \la_i)$ divides $\minp_a$, and $\chpoly_a$ divides $\minp_a$.

Since $\minp_a$ divides $\chpoly_a$ by Cayley-Hamilton, the corollary is proved.
\end{proof}

\begin{prop} \label{alt.sum}
For $a$ and $a'$ in $M_n(F)$, the following are true:
   \begin{enumerate}
  \item $a$ is invertible if and only if $c_n(a) \ne 0$.
  \item {\rm (Jacobi formula)} $c_n(a) = \sum_{\s \in S_n} (\sign \s) a_{1\s(1)} a_{2\s(2)} \cdots a_{n\s(n)}$.
  \item $c_n(aa') = c_n(a) c_n(a').$
   \end{enumerate}
\end{prop}

\begin{proof} 
(1): The matrix $a$ is invertible if and only if the kernel of $a$ is trivial, if and only if $0$ is not an eigenvalue of $A$.   By Theorem \ref{same.thm}, this is true if and only if $c_n(a) \ne 0$.

\smallskip

We now follow \cite[\S9]{Axler:dwd}.  Write $d(a)$ for the right-hand side of (2).  A straightforward rearrangement of terms as in \cite[p.~179]{BW} 
shows that $d(aa') = d(a) d(a')$.  Therefore, (2) implies (3).

\smallskip

We now prove (2).  Suppose first that $a$ is upper-triangular.  Then both sides of (2) are just the product of the diagonal entries of $a$, hence (2) holds in this case.

Now consider the general case.  Since $c_n(a)$ and $d(a)$ are unchanged if we enlarge our base field, we may assume that $F$ is algebraically closed and hence that $a$ is similar to an upper triangular matrix $a'$, i.e., $a' = bab^{-1}$ for some $b \in M_n(F)$.  Then
   \[
c_n(a') = d(a') = d(bab^{-1}) = d(b) d(ab^{-1}) = d(ab^{-1})d(b) = d(ab^{-1} b) = d(a).
  \]
Since $c_n(a) = c_n(a')$ by \ref{Jac.prop}.2, we have proved (2).
\end{proof}

\begin{defn} \label{det.defn}
For $a$ an element of a finite-dimensional $F$-algebra $A$, we define the {\em trace of $a$} to be
   \[
  \tr_{A/F}(a) := c_1(a)
  \]
and the {\em determinant of $a$} to be
  \[
  \det_{A/F}(a) := c_n(a).
  \]
If there is no danger of ambiguity, we write simply $\det_A$ or $\det$ instead of $\det_{A/F}$ and similarly for the trace.
\end{defn}

By Theorem \ref{same.thm}, the trace $\tr_{M_n(F)}$ and determinant $\det_{M_n(F)}$ are the usual trace and determinant from linear algebra.

\section{Quaternions}

\begin{eg} \label{quats.eg}
Hamilton's quaternions --- usually denoted by $\H$ --- are defined to be the ring constructed by taking the complex numbers $\C$ and adjoining an element $j$ such that $j^2 = -1$, $j$ commutes with real numbers, and $ij = -ji$.  Note that $\H$ is an $\R$-algebra but not a $\C$-algebra, since $\C$ is not in the center of $\H$.  It has $\R$-basis 1, $i$, $j$, $k$, where $k = ij$.  A lot of interesting information about the quaternions can be found in \cite[Ch.~7]{Numbers}.

Set $\phi$ to be the $\R$-linear map $\H \ra M_2(\C)$ defined by
    \begin{gather*}
    1 \mapsto \tbtmat{1}{0}{0}{1}  \text{for $r$ real}, \quad i \mapsto \tbtmat{i}{0}{0}{-i},\\
     \quad j \mapsto \tbtmat{0}{1}{-1}{0},  \quad \text{and}  \quad k \mapsto \tbtmat{0}{i}{i}{0}.
  \end{gather*}
This extends to an isomorphism $\phi \!: \H \ot \C \ra M_2(\C)$ (as $\C$-algebras).  


Every quaternion  can be written as $q = r + si + uj + vk = z + wj$ for some real numbers $r$, $s$, $u$, $v$ and complex numbers $z = r + si$ and $w = u+ vi$.
We have
   \[
   \phi(q) = \tbtmat{z}{w}{-\overline{w}}{\overline{z}}.
  \]
Since the characteristic polynomial is unchanged when we enlarge our base field, we find that
  \[
  \tr_\H(q) = \tr_{M_2(\C)}(\phi(q)) = z + \overline{z} = 2r
  \]
and
  \[
  \det_\H(q) = \det_{M_2(\C)}(\phi(q)) = z \overline{z} + w \overline{w} = r^2 + s^2 + u^2 + v^2.
  \]
\end{eg}

\begin{eg}[Matrices over the quaternions] \label{Cay.eg}
Write $M_2(\H)$ for the set of 2-by-2 matrices with entries in $\H$.  The obvious addition and multiplication make it into a 16-dimensional $\R$-algebra.

In \cite{Cayley}, Cayley defined a determinant 
  \[
  \cdet  \!: M_2(\H) \ra \H \quad \text{by} \quad \cdet \tbtmat{q_{11}}{q_{12}}{q_{21}}{q_{22}} =
   q_{11} q_{22} - q_{21} q_{22}
   \]
He noted that his determinant has some unsavory properties, for example that
   \[
   \cdet \tbtmat{q}{q'}{q}{q'} = 0 \quad \text{for all $q$, $q' \in \H$,}
   \]
but
  \[
  \cdet \tbtmat{i}{i}{j}{j} = ij - ji = 2ij \ne 0.
  \]
  
Let us constrast this with the determinant that we have just defined.
Just as for $\H$ above, there is an isomorphism $\phi_2 \!: M_2(\H) \ot \C  \ra M_4(\C) $ such that
   \[
  \phi_2 \tbtmat{q_{11}}{q_{12}}{q_{21}}{q_{22}} = \tbtmat{\phi(q_{11})}{\phi(q_{12})}{\phi(q_{21})}{\phi(q_{22})}.
  \]
(Recall that $\phi(q)$ is a 2-by-2 complex matrices for every $q \in \H$.) 
If $m \in M_2(\H)$ has a repeated row or column, then so does $\phi_2(m)$, hence
   \[
  \det_{M_2(\H)}(m) = \det_{M_4(\C)}(\phi_2(m)) = 0.
  \]
More generally, our trace and determinant have the nice properties of the usual trace and determinant for matrices as given in \ref{Jac.prop}.

For a more comprehensive discussion of various types of determinants for $M_2(\H)$, see \cite{Aslaksen}.  Aslaksen refers to our $\det_{M_2(\H)}$ as the {\em Study determinant}.
\end{eg}

\begin{eg}[Central simple algebras]
A typical topic for a first year graduate algebra
course is Wedderburn's description of simple artinian rings: they
are isomorphic to $M_r(D)$ for $D$ a skew field.  Write $F$ for
the center of $D$ (which is necessarily a field), and suppose that
$D$ is finite-dimensional over $F$.  Such an algebra $M_r(D)$ is called {\em central simple}.
We have just seen two examples of these, with $F = \R$, $D = \H$, and $r = 1$, 2.

The trace $\tr_{M_r(D)}$
and determinant $\det_{M_r(D)}$ are called the {\em
reduced trace} and {\em reduced norm} respectively.   They are 
usually constructed by ``Galois descent'' as in
\cite[p.~145]{Draxl} or \cite{Pierce}, but here we get them as a
consequence of the existence of the characteristic polynomial, which is true for a much broader class of algebras.
\end{eg}

\section{More properties of the characteristic polynomial}

Here we discuss the example of finite-dimensional field extensions (\ref{fields.eg}) and prove some more nice properties of the characteristic polynomial (\ref{Jac.cor}).  

Write $\EndF(A)$ for the set of $F$-linear maps $A \ra A$.  It is an $F$-algebra; its multiplication is function composition.  It is isomorphic to $M_m(F)$.

For $a \in A$, write $L_a$ for the element of $\EndF(A)$ defined by
   \[
  L_a(b) = ab  \quad \text{for $b \in A$.}
  \]
The map $a \mapsto L_a$ defines an $F$-algebra homomorphism called the {\em left regular
representation} of $A$.  This homomorphism is injective: if $L_a = 0$, then $L_a(a') = 0$ for all $a' \in A$, hence $0 = L_a(1_A) = a \cdot 1_A = a$.

\begin{eg} \label{bbk.eg}
Let $A$ be $M_n(F)$ or more generally a simple ring with center
$F$ such that $\dim_F A = n^2$.  One can show that
   \begin{equation} \label{power.reln}
   \gminp_{L_a, \EndF(A)/F} = (\gminp_{a,A/F})^n
   \end{equation}
for every $a \in A$, see e.g.~\cite[Ch.~VIII, \S12.3]{Bbk:alg}.
\end{eg}

\begin{prop} \label{factor} \label{divides.2}
Let $a$ be an element in a finite-dimensional $F$-algebra $A$.
The minimal polynomial $\minp_{a/F}$ divides the characteristic polynomial $\gminp_{a,A/F}$ of $a$ which divides the characteristic polynomial $\gminp_{L_a,\EndF(A)/F}$, all in $F[x]$.
All three polynomials have the same irreducible factors in $F[x]$.
\end{prop}

\begin{proof}
Since $\gminp_{a}(a) = 0$ by Cayley-Hamilton (\ref{Jac.prop}.1), the minimal polynomial of $a$ divides the characteristic polynomial $\gminp_a$.  Since the left regular representation is injective, we have that $\gminp_a$ divides $\gminp_{L_a}$ by \ref{Jac.prop}.5.

We are reduced to showing that $\gminp_{L_a}$ and $\minp_a$ have the same irreducible factors.  Since the left regular representation is injective, $a$ and $L_a$ have the same minimal polynomials.  That is, we need only show that $\gminp_{L_a}$ and $\minp_{L_a}$ have the same irreducible factors.  Since $\EndF(A)$ is isomorphic to $M_m(F)$ for $m = \dim_F A$, we are done by \ref{divides.1}.
\end{proof}

The proposition gives us the power to handle another example.

\begin{eg}[Finite-degree field extensions] \label{fields.eg}
Let $A$ be an extension field of $F$ of finite dimension $m$.  Every element $a \in A$ gives an element $L_a$ in $\EndF(A) \cong M_m(F)$, and the characteristic polynomial of $a$ divides the characteristic polynomial of $L_a$ by the proposition.  The trace and norm of $a$ are typically defined to be the trace and determinant of $L_a$.

\smallskip

If $A$ is separable over $F$, then by the Theorem of the Primitive Element, $A = F[\theta]$ for some $\theta \in A$.  The minimal polynomial of $\theta$ has degree $m$.  Since it divides the characteristic polynomial of $L_\theta$ by \ref{divides.2} and that polynomial has degree $m$, we find that 
   \[
   \chpoly_{\theta,A/F} = \chpoly_{L_\theta}.
   \]
Since the characteristic polynomial has the same degree for all $a \in A$, we have
   \[
   \chpoly_a = \chpoly_{L_a} \quad \text{for all $a \in A$.}
   \]
In particular, for finite separable field extensions, our trace and determinant agree with the usual trace and norm.

\smallskip

If $A$ is not separable over $F$, there can be some disagreement.  For example, let $F = \F_2(u,v)$ where $\F_2$ is the field with 2 elements and $u$, $v$ are indeterminates.  The field $A = F(\sqrt{u}, \sqrt{v})$ is a purely inseparable extension of degree 4 with $F$-basis $1$, $\sqrt{u}$, $\sqrt{v}$, $\sqrt{uv}$.  The generic element
   \[
  \gamma = t_1 \cdot 1 + t_2  \sqrt{u} + t_3  \sqrt{v} + t_4  \sqrt{uv}
   \]
has minimal polynomial
  \[
  \gamma^2 - (t_1^2 + t_2^2 u + t_3^2 v + t_4^2 uv).
  \]
Here the characteristic polynomial of each element $a \in A$ divides but is not equal to the characteristic polynomial of $L_a$.
\end{eg}

The proposition also
allows us to prove that many nice properties of the characteristic
polynomial of a matrix also hold for characteristic polynomials of elements of $A$.

Recall (\ref{deg.rmk}) that the characteristic polynomial has the same degree for every $a \in A$.  We say that $A$ {\em has degree $n$} if $\chpoly_a$ has degree $n$.

\begin{cor} \label{Jac.cor}
Let $A$ be a finite-dimensional $F$-algebra of degree $n$.  For $a \in A$, we have:
   \begin{enumerate}
   \item $c_i(1_A) = \binom{n}{i}$.  In particular, $\tr_A(1_A) = n$ and $\det_A(1_A) = 1$.
   \item $a$ is invertible if and only if $\det_A(a) \ne 0$.
   \item $a$ is nilpotent if and only if $\chpoly_a = x^n$.
   \end{enumerate}
\end{cor}

\begin{proof}
(1): The minimal polynomial of $1_A$ is $x - 1$, hence Proposition \ref{factor} gives that
   \[
   \gminp_{1_A} = (x - 1)^n.
   \]

\smallskip

(2, $\Leftarrow$): Define the {\em adjoint} of $a$, denoted by $\adj a$, to be
   \[
   \adj a = (-1)^{n + 1} [ a^{n-1} - c_1(a) a^{n-2} + \cdots + (-1)^{n-1} c_{n-1}(a) ].
  \]
Then
   \[
  a \cdot \adj a = (-1)^{n+1}[ \gminp_a(a) - (-1)^n c_n(a)] = \det_A(a) 1_A.
   \]
hence, if $\det_A(a)$ is not zero, $a$ is invertible.

\smallskip

(2, $\Rightarrow$): If $a$ is invertible, then $L_a$ is invertible with inverse $L_{a^{-1}}$.  By \ref{alt.sum}.1, the constant term $\det_{\EndF(A)}(L_a)$ of $\gminp_{L_a}$ is not zero.  Since $\gminp_{a,A/F}$ and $\gminp_{L_a,\EndF(A)/F}$ have the same irreducible factors in $F[x]$, the constant term $\det_A(a)$ of $\gminp_{a,A/F}$ is not zero.

\smallskip

(3): $a$ is nilpotent if and only if it satisfies $x^r$ for some natural number $r$, if and only if it has minimal polynomial $x^p$ for some natural number $p$.  Since the minimal polynomial and characteristic polynomial have the same irreducible factors \ref{factor}, this holds if and only if the characteristic polynomial of $a$ is $x^n$.
\end{proof}


\begin{rmk} \label{Sch}
One might be tempted to accept the traditional definition of characteristic polynomial for matrices, and then define the characteristic polynomial of $a \in A$ as $\chpoly_{L_a, \EndF(A)/F}$.  (This is logically equivalent to the usual definition of the norm and trace  in the case where $A$ is a finite-degree field extension of $F$, see e.g.\ \cite[VI.5.6]{Lang}.)  However, Example \ref{bbk.eg} shows that in some cases one wants to take an $n$-th root of $\chpoly_{L_a}$.  But there is a more serious problem: There is no strong mathematical reason to prefer the left regular representation over the right regular representation (defined in the obvious
 manner as $a \mapsto R_a$), and the characteristic polynomials of
 $L_a$ and $R_a$ may differ.  Adrian Wadsworth points out
 that the upper triangular matrices from Example \ref{ut.eg}
 provide an example of this difficulty.  In particular,
 the generic element $\gamma$ has
 \[
 \gminp_{L_a}(\gamma) = (x-t_1)^n (x-t_2)^{n-1} \cdots (x-t_{n-1})^2
 (x-t_n)
 \]
and
 \[
 \gminp_{R_a}(\gamma) = (x-t_1) (x-t_2)^2 \cdots (x-t_{n-1})^{n-1}
 (x-t_n)^n.
 \]
\end{rmk}

\section{The product formula for determinants} \label{prod.sec}

In this section, we prove that the usual product formula for determinants of matrices holds for an arbitrary finite-dimensional $F$-algebra $A$:

\begin{thm}  \label{product}
For every $a$, $a' \in A$, we have
   \[
  \det_A(aa') = \det_A(a)\,  \det_A(a').
  \]
\end{thm}

We postpone the proof until the end of the section.  A somewhat different proof can be found in \cite[\S{VI.5}]{Jac:J}, but it is even more complicated and arcane than what you will find below.

\borelless \label{poly.fcn}
Recall that $\det_A$ is both a function $A \ra F$ and an element of $R = F[t_1, \ldots, t_m]$.  To evaluate the function $\det_A$ on $a \in A$, we write $a$ in terms of our basis $a_1$, $\ldots$, $a_m$ as $a = \sum_i \alpha_i a_i$ and substitute $t_i \mapsto \alpha_i$ in the polynomial $\det_A$.  In this manner, we may view each polynomial in $R$ as a function $A \ra F$.

\begin{lem} \label{product.lem}
Suppose that $F$ is infinite.  Let $f \in R$ be such that 
   \[
   f(1_A) = 1 \quad \text{and} \quad \text{$f(aa') = f(a)\,f(a')$ for all $a, a' \in A$.}
   \]
If $g \in R$  divides  $f$ and has $g(1_A) = 1$, then
   \[
   g(aa') = g(a) \,g(a') \quad \text{for all $a$, $a' \in A$.}
   \]
\end{lem}

\begin{proof}
Since $f(1_A) = g(1_A) = 1$, we have $f,  g \ne 0$ in $R$.  Suppose that $f$ or $g$ is a unit in $R$, i.e., is in $R^* = F^*$.  Then $g$ is a unit, hence $g = 1$ and the lemma holds.  

We may assume that $f$ and $g$ are nonzero nonunits.  Since $R$ is a UFD, we may write $f = f_1 f_2 \cdots f_r$ where $f_k$ is irreducible in $R$ for every $k = 1$, $\ldots$, $r$.  By multiplying $f_k$ by an element of $F^*$ if necessary, we may assume that $f_k(1_A) = 1$ for every $k$.  Since $g$ divides $f$ and $g(1_A) = 1$, $g$ is a product
   \[
  g = \prod_{\ell \in L} f_\ell  \quad \text{for some $L \subseteq [1, n]$.}
  \]
To prove the lemma, it suffices to prove that $f_k(aa') = f_k(a) f_k(a')$ for all $k \in [1,n]$ and $a$, $a' \in A$.

Set $R' = F[u_1, \ldots, u_m, v_1, \ldots, v_m]$ for $u_1, \ldots, u_m$, $v_1, \ldots, v_m$ independent indeterminates, and let $F'$ be the quotient field of $R'$.  Set
  \[
  \mu = \sum_i u_i a_i \quad \text{and} \quad \nu = \sum_j v_j a_j \quad \text{in $A \ot R'$.}
  \]
  
As in \ref{poly.fcn}, every element of $R$ defines a map $A \ot R' \ra R'$.  For example, to find $f(\mu)$, one takes $f$ and substitutes $t_i \mapsto u_i$.  Similarly, one gets $f(\nu)$ by substituting $t_j \mapsto v_j$.

We claim that $f(\mu \nu) = f(\mu) f(\nu)$.  Let $\delta := f(\mu \nu) - f(\mu) f(\nu)$ in $R'$.  Every element of $R'$ defines a function $F^{2m} \ra F$ by plugging in for the the $u$'s and $v$'s.  Substituting $u_i \mapsto \alpha_i$ and $v_j \mapsto \beta_j$ into $\delta$ for $\alpha_i$, $\beta_j \in F$, we obtain
   \[
  f(ab) - f(a) f(b) \quad \text{for $a = \sum \alpha_i a_i$ and $b = \sum \beta_j a_j$,}
  \]
which is 0 by hypothesis.  That is, $\delta$ gives the map $F^{2m} \ra F$ which is identically 0.  Since $F$ is infinite, $\delta = 0$, which proves the claim.

Thus 
   \[
   \prod_{k=1}^r f_k(\mu \nu) = f(\mu \nu) = f(\mu) f(\nu) = \prod_{k=1}^r f_i(\mu) f_i(\nu) \quad \text{in $R'$.}
   \]
Recall that one obtains $f_k(\mu)$ and $f_k(\nu)$ by substituting one set of indeterminates for another in $f_k$.  Hence, since $f_k$ is irreducible in $R$, the polynomials $f_k(\mu)$, $f_k(\nu)$ are irreducible in $R'$ for all $k$.  Therefore the prime factorization of $f_k(\mu \nu)$ in the UFD $R'$ is a product of $f_p(\mu)$ and $f_q(\nu)$ for some $p$'s and $q$'s.  Substituting in for the $v_j$'s so that $\nu$ is sent to $1_A$ maps
  \[
  f_k(\mu \nu) \mapsto f_k(\mu), \quad f_p(\mu) \mapsto f_p(\mu), \quad \text{and} \quad f_q(\nu) \mapsto f_q(1_A) = 1
   \]
for all $p$ and $q$.  Hence the only irreducible factor of $f_k(\mu \nu)$ amongst the $f_p(\mu)$ terms is equal to $f_k(\mu)$.  Similarly, substituting in for the $u_i$'s so that $\mu \mapsto 1_A$, we have $f_k(\mu \nu) \mapsto f_k(\nu)$, hence
   \[
  f_k(\mu \nu) = f_k(\mu) f_k(\nu) \quad \text{for $k = 1, \ldots, n$.}
  \]
By substituting $\mu \mapsto a$ and $\nu \mapsto a'$, we obtain that $f_k(aa') = f_k(a) f_k(a')$ for all $a$, $a' \in A$.
\end{proof}

\begin{proof}[Proof of Theorem \ref{product}]
Since the determinant of an element of $A$ is unchanged when we enlarge the base field (\ref{scalar.ext}), we may assume that $F$ is infinite.
Let $\gamma = \sum t_i a_i$ be a generic element of $A$ as in \S\ref{chpoly.sec}, and consider the element $L_\gamma$ in $\EndK(A \ot K)$.

The minimal polynomial $\minp_\gamma$ divides $\gminp_{L_\gamma,\EndK(A \ot K)/K}$ in $K[x]$ by \ref{factor}.  By \ref{poly.lem}, both polynomials actually lie in $R[x]$ for $R = F[t_1, \ldots, t_m]$, hence $\minp_\gamma$ divides $\gminp_{L_\gamma}$ in $R[x]$ by Gauss' Lemma.  Substituting $x \mapsto 0$ defines a surjection $R[x] \ra R$ which sends $\minp_\gamma \mapsto (-1)^d \det_A$ for $d = \deg(\minp_\gamma)$ and $\gminp_{L_\gamma}$ to the function $a \mapsto (-1)^m \det_{\EndF(A)}(L_a)$.  Since substitution is a ring homomorphism, $\det_A$ divides $a \mapsto \det_{\EndF(A)}(L_a)$ in $R$.

We have 
  \[
  \det_{\EndF(A)}(L_{aa'}) = \det_{\EndF(A)}(L_a L_{a'}) = \det_{\EndF(A)}(L_a) \ \det_{\EndF(A)}(L_{a'}),
  \]
where the last equality is by \ref{alt.sum}.3 since $\EndF(A) \cong M_m(F)$.  Note that   
  \[
  \det_{\EndF(A)}(L_{1_A}) = \det_{M_m(F)}(1_{M_m(F)}) = 1
  \]
   and $\det_A(1_A) = 1$ by \ref{Jac.cor}.1.  Since $F$ is infinite, Lemma \ref{product.lem} gives:
   \[
  \det_A(aa') = \det_A(a) \  \det_A(a')
  \]
as desired.
\end{proof}

To summarize, we defined the characteristic polynomial of an element in a finite-dimensional $F$-algebra in \S\ref{chpoly.sec}.  We defined the determinant $\det_A$ to be the constant term of this polynomial (\ref{det.defn}).  In the case $A = M_n(F)$, we found that $\det_{M_n(F)}$ is given by the Jacobi formula (\ref{alt.sum}.2), hence the product formula holds for $\det_{M_n(F)}$ (\ref{alt.sum}.3).  Finally, we used Prop.~\ref{divides.2} to prove the product formula for $\det_A$ (\ref{product}). 

\section{Miscellaneous remarks} \label{misc}

This section is a survey of other results.  It is necessarily briefer and more technical than the rest of the paper.

\mysubsection{The ``usual'' definition of the characteristic polynomial} 
For $a \in A$, I claim that the formula 
   \begin{equation} \label{ch.usdef}
     \chpoly_a = \det(x \cdot 1_A - a) 
  \end{equation}
holds in $F[x]$.  But what does the expression ``$\det(x \cdot 1_A - a)$'' mean?  For matrices, the determinant is given by the Jacobi formula \ref{alt.sum}.2 involving the entries of a matrix; this formula makes sense whether the entries are elements of $F$ or polynomials in $F[x]$.  The same reasoning holds for our more general notion of determinant.

One way to prove \eqref{ch.usdef} is as follows.  First prove it for $A$ the algebra of upper triangular matrices from Example \ref{ut.eg}.  Then prove \eqref{ch.usdef} for $A = M_n(F)$ by reducing to the upper triangular case as in the proof of Th.~\ref{same.thm}.  This implies \eqref{ch.usdef} for general $A$ by the arguments in \cite[p.~225]{Jac:J}.

\mysubsection{$\mathbf{ch.poly.}_{ab} = \mathbf{ch.poly.}_{ba}$} If $a$ and $b$ are $n$-by-$n$ matrices, it is well-known that $\chpoly_{ab} = \chpoly_{ba}$.  We can prove it directly for general $A$.  Let $c_i$ be one of the coefficients of the characteristic polynomial as in \eqref{char.poly}, so $c_i \!: A \ra F$ and we want to show that $c_i(ab) = c_i(ba)$ for all $i$.

First assume that $b$ is invertible.  Then 
   \[
   c_i(ab) = c_i(b(ab)b^{-1}) = c_i(ba)
   \]
by \ref{Jac.prop}.2.  

Now the set $U$ of pairs $(a, b) \in A \times A$ such that $b$ is invertible is a nonempty open subset of $A \times A$ in the Zariski topology.  The map $A \times A \ra F$ defined by
   \[
  (a, b) \mapsto c_i(ab) - c_i(ba)
  \]
is given by a polynomial in the coordinates of $a$ and $b$ and is zero on $U$, hence is identically zero on $A \times A$.  This proves that $c_i(ab) = c_i(ba)$ for all $a, b \in A$ and all $i$.

\mysubsection{The Pfaffian}
Suppose that $\chr F \ne 2$.  Write $\Skew_n(F)$ for the vector space of {\em skew-symmetric} $n$-by-$n$ matrices, i.e., matrices $a$ such that $a^t = -a$.  

If $n$ is odd, then $\det(-a) = -\det(a)$ by \ref{Jac.prop}.3.  Since
$\det(a^t) = \det(a)$ by \ref{Jac.prop}.2, all $n$-by-$n$ skew-symmetric matrices have determinant 0.

Suppose that $n$ is even.  There is a polynomial map
  \begin{equation} \label{Pf.eq}
  \Pf \!: \Skew_n(F) \ra F \quad \text{such that} \quad (\Pf(a))^2 = \det_{M_n(F)}(a)
  \end{equation}
called the {\em Pfaffian},  see e.g.\ \cite[\S{XV.9}]{Lang}.  It seems a bit mysterious!

(Note that Equation \eqref{Pf.eq} only determines the Pfaffian up to sign; classically one chooses an invertible skew-symmetric matrix $S$ with $\det(S)  = 1$ and fixes the sign of $\Pf$ so that $\Pf(S) = 1$.)

In fact, the Pfaffian exists as a consequence of the characteristic polynomial as defined in \S\ref{chpoly.sec}.  Fix a matrix $S$ as in the preceding paragraph and define a multiplication $\cdot$ on $\Skew_n(F)$ given by
     \[
     a \cdot b = \frac12 (a S^{-1} b + b S^{-1} a).
     \]
This makes $\Skew_n(F)$ into a Jordan $F$-algebra with identity element $S$.  These algebras arise naturally in the classification of central simple Jordan algebras\footnote{See \cite[\S{V.7}]{Jac:J}.  They correspond to Lie algebras of type C in the Killing-Cartan classification.}; they are associated with a nondegenerate skew-symmetric bilinear form on $F^n$.  

As mentioned in the introduction, the theory developed in \S\ref{chpoly.sec} can be extended to give a characteristic polynomial and a determinant for this Jordan algebra.  One finds that the determinant {\em is} the Pfaffian:
   \[
     \det_{\Skew_n(F)}(a) = \Pf(a),
    \]
see \cite[pp.~230--232]{Jac:J}.  

\mysubsection{Matrices over rings} Apparently, we have given a canonical construction of an analogue
of the characteristic polynomial and hence the determinant for any
algebra over a field $F$.  But
the characteristic polynomial and determinant are typically defined for
$M_n(R)$ where $R$ is merely a commutative ring with 1; our definition also works in that case with some modification.

First, consider $M_n(\Z)$ as a $\Z$-algebra, i.e., the case $A = M_n(\Z)$ and $F = \Z$.  Since $\Z$ is a noetherian UFD and $M_n(\Z)$ is a free $\Z$-module of finite rank, the results in \S\ref{chpoly.sec} hold with this $A$ and $F$.  

Lemma \ref{scalar.ext} holds if we require $E$ to be a noetherian UFD containing $F$
so that we may define the characteristic polynomial over $E$.  (Since
$A$ is a free $F$-module, $A \otimes E$ is a free $E$-module.) Let
$L$ be the field of rational functions in $t_1$, $\ldots$, $t_m$ over the quotient field of $E$.  The minimal polynomial for $\gamma$ over $K$ is
also the minimal polynomial for $\gamma$ over $L$ by
\ref{scalar.ext}, hence characteristic polynomials are also
unchanged by scalar extension in this more general setting.  In particular, the characteristic polynomial of a matrix $a \in M_n(\Z)$ is the same as the characteristic polynomial of $a$ considered as a matrix in $M_n(\Q)$, which is the usual characteristic polynomial by Th.~\ref{same.thm}.

A coefficient $c_i$ of the characteristic polynomial as in \eqref{char.poly} is given by a polynomial in $\Z[t_1, \ldots, t_m]$, and it is evaluated on a matrix $a$ as described in \ref{poly.fcn}.  For example, undergraduates are taught how to evaluate $c_n = \det$, which is given by the Jacobi formula \ref{alt.sum}.2 with respect to the standard basis of $M_n(\Z)$.  Just as for the determinant, one can evaluate $c_i$ on a matrix in $M_n(\Z) \ot R = M_n(R)$ for $R$ any commutative ring with 1, and the resulting characteristic polynomial has the usual felicitous properties.


\providecommand{\bysame}{\leavevmode\hbox to3em{\hrulefill}\thinspace}

\end{document}